\documentclass[english,twoside,11pt]{article}

\usepackage[german,english]{babel}
\usepackage{amsmath}
\usepackage{amssymb,latexsym}
\usepackage{url}
\usepackage{graphicx}
\usepackage{booktabs}
\usepackage{psfrag}
\usepackage{paralist}

\newtheorem{deff}{Definition}

\newtheorem{prop}[deff]{Proposition}
\newtheorem{thm}[deff]{Theorem}
\newtheorem{prob}[deff]{Problem}
\newtheorem{cor}[deff]{Corollary}
\newtheorem{conj}[deff]{Conjecture}

\makeatletter

\title{Isomorphism-Free Lexicographic Enumeration of Triangulated Surfaces and $3$-Manifolds}

\author{\Large Thom~Sulanke and Frank H.~Lutz\footnote{Supported by the DFG Research Group ``Polyhedral Surfaces'', Berlin}}

\date{}

\begin{document}

\selectlanguage{english}

\maketitle

\begin{abstract}
We present a fast enumeration algorithm for combinatorial $2$- and $3$-manifolds. 
In particular, we enumerate all triangulated surfaces with $11$ and $12$ vertices 
and all triangulated $3$-manifolds with $11$ vertices. We further
determine all equivelar polyhedral maps on the non-orientable surface of genus $4$
as well as all equivelar triangulations of the orientable surface of
genus $3$ and the non-orientable surfaces of genus $5$ and $6$.
\end{abstract}

\section{Introduction}

Triangulations of manifolds with \emph{few} vertices provide a valuable source 
of interesting and extremal testing examples for various conjectures and 
open problems in combinatorics, optimization, geometry, and topology:
\begin{compactitem}
\item \emph{How many vertices are needed to triangulate a given manifold?}
\item \emph{What do the face vectors of simplicial spheres and manifolds look like?}
\item \emph{Is a given triangulation of a sphere shellable, polytopal,
            does it satisfy the Hirsch conjecture?} 
\item \emph{Is a given triangulated surface geometrically realizable
            as a polyhedron in $3$-space?}
\item \emph{Do simplicial $2$-spheres have polytopal realizations with
            small coordinates?}
\end{compactitem}

According to Rado \cite{Rado1925} and Moise \cite{Moise1952},
(closed, compact) $2$- and $3$-manifolds can always be triangulated
as (finite) simplicial complexes. Moreover, triangulated $2$- and $3$-manifolds
always are \emph{combinatorial manifolds}, i.e., triangulated
manifolds such that the links of all vertices are standard PL (i.e.,
piecewise linear) spheres. 
It immediately follows that triangulations of $2$- and $3$-manifolds 
\emph{can be enumerated}: For any given positive integer~$n$ 
we can produce in a finite amount of time a complete list (up to combinatorial
isomorphism) of all triangulated $2$- respectively $3$-manifolds with
$n$ vertices. 

For example, a conceptually simple, but highly inefficient 
enumeration approach would be to first generate all $2^{\binom{n}{3}}$
pure $2$-dimensional (respectively all $2^{\binom{n}{4}}$
pure $3$-dimensional) simplicial complexes. These complexes
would then be tested in a second step to determine whether or not all links 
are triangulated circles (triangulated $2$-spheres), which can be
done purely combinatorially. In a third step, isomorphic copies of triangulations,
i.e., copies that can be transformed into each other by relabeling
the vertices, would be identified.

It is the second step which fails to work in higher dimensions:
There is \emph{no} algorithm to decide whether a given 
$r$-dimensional simplicial complex is a PL $r$-sphere if $r\geq 5$;
cf.~\cite{VolodinKuznetsovFomenko1974}.
For $r=4$ it is unknown whether there are algorithms to recognize PL $4$-spheres.
For $r=3$ there \emph{are} algorithms to recognize the $3$-sphere 
(see \cite{Rubinstein1995}, \cite{Thompson1994}, and also
\cite{King2004}, \cite{Matveev1995}, \cite{Mijatovic2003}),
however, all the known algorithms are exponential and hopeless to implement.
Therefore, in principle, combinatorial $4$-manifolds (with $3$-dimensional links)
can be enumerated, whereas the enumeration problem for combinatorial $5$-manifolds
is open, and there is no enumeration algorithm for combinatorial
$(r+1)$-manifolds for $r\geq 5$.

We are, of course, not only interested in the combinatorial types
of triangulated manifolds: In an additional step we want to
determine the topological types of the triangulations
obtained by the enumeration. 
Algorithmically, it is easy to figure out the topological type 
of a triangulated surface (by computing its Euler characteristic 
and its orientability character). As mentioned before, 
there are algorithms to recognize the $3$-sphere, 
and it is even possible to recognize Seifert manifolds \cite{Mijatovic2004}. 
For general $3$-manifolds, however, there are no algorithmic tools 
available yet to determine their topological types (although 
Perelman's proof \cite{Perelman2003bpre} of Thurston's 
geometrization conjecture \cite{Thurston1982}  
gives a complete classification of the geometric types of $3$-manifolds). 
In particular, hyperbolic $3$-manifolds are difficult to deal with. 
For triangulations with few vertices it turned out
that heuristics (e.g., \cite{BjoernerLutz2000}, \cite{Matveev2005}) 
can be used for the recognition, thus allowing for a complete 
classification of the topological types of the examples obtained
by the enumeration.

At present, there are three major enumeration approaches known
to generate triangulated manifolds (see the overview \cite{Lutz2005apre}): 
\begin{compactitem}
\item \emph{generation from irreducible triangulations} 
      (\cite{BrinkmannMcKay2007}, \cite{Sulanke2006apre}, \cite{Sulanke2006bpre},
      with the programs \texttt{plantri} \cite{plantri} 
      of Brinkmann and McKay and \texttt{surftri}
      \cite{Sulanke2005pre} of Sulanke
      implementing this approach),
\item \emph{strongly connected enumeration}
      (\cite{LutzSullivan2005pre};
      \cite{AltshulerBokowskiSchuchert1996}, \cite{Bokowski2006pre}),
\item and \emph{lexicographic enumeration}
      (\cite{EllinghamStephens2005}, \cite{Lutz2005apre}, \cite{Lutz2006apre}; 
      \cite{KoehlerLutz2005pre}, \cite{KuehnelLassmann1985-di}).
\end{compactitem}

A triangulation is \emph{irreducible} if it has no \emph{contractible edge},
i.e., if the contraction of any edge of the triangulations yields a
simplicial complex, which is not homeomorphic to the original
triangulation. According to Barnette and Edelson \cite{BarnetteEdelson1988},
every surface has only finitely many irreducible triangulations from
which all other triangulations of the surface can be obtained by a suitable
sequence of vertex splits. In this manner, triangulations of a
particular surface with $n$ vertices can be obtained in two steps
by first generating all irreducible triangulations of the surface 
with up to $n$ vertices, from which further triangulations with $n$
vertices are obtained fast by vertex splits; see \cite{BrinkmannMcKay2007}, 
\cite{Sulanke2006apre}, \cite{Sulanke2006bpre}. Unfortunately, every 
$3$-manifold has infinitely many irreducible triangulations; cf.\
\cite{DeyEdelsbrunnerGuha1999}. Even for surfaces, the generation
of the finitely many irreducible triangulations is difficult,
with complete lists available only for the $2$-sphere, the $2$-torus,
the orientable surface of genus $2$, 
and the non-orientable surfaces of genus up to $4$; see \cite{Sulanke2006apre}, \cite{Sulanke2006bpre}
and the references contained therein.

Strongly connected enumeration, in particular, turned out to be successful
for the enumeration of triangulated $3$-manifolds with small edge degree \cite{LutzSullivan2005pre},
but is otherwise not very systematic.

The third approach, lexicographic enumeration, generates triangulations in \emph{canonical form},
that is, for every fixed number $n$ of vertices a lexicographically sorted list of 
triangulated manifolds is produced such that every listed triangulation 
with $n$ vertices is the lexicographically smallest set of triangles (tetrahedra) 
combinatorially equivalent to this triangulation and is lexicographically smaller 
than the next manifold in the list. 

In this paper, we present an improved version of the algorithm 
for lexicographic enumeration from \cite{Lutz2005apre}.
The triangulations are now generated in an isomorphism-free way; 
see the next section for a detailed discussion. This improvement 
led to a substantial speed up of the enumeration. In particular,
with the implementation \texttt{lextri} of the first author,
we were able to enumerate all triangulated surfaces
with $11$ and $12$ vertices (Section~\ref{sec:2d}) and
all triangulated $3$-manifolds with $11$ vertices (Section~\ref{sec:3d}).
Moreover, we enumerated all equivelar triangulations 
of the orientable surface of genus $3$ and of 
the non-orientable surfaces of genus $4$, $5$, and $6$ (Section~\ref{sec:equiv}).

\section{Isomorphism-Free Enumeration}
\label{sec:free}

It is a standard problem with algorithms for the enumeration of particular combinatorial objects
to avoid isomorphic copies of the objects as early as possible during their generation; 
see Read \cite{Read1978} and McKay \cite{McKay1998} for a general discussion.

Our aim here is to give an isomorphism-free enumeration algorithm 
for triangulated surfaces and $3$-manifolds with a fixed number $n$ of
vertices $1,2,\dots,n$. The algorithm is based on lexicographic enumeration as discussed in \cite{Lutz2005apre}.
For simplicity, we describe the algorithm for surfaces, however,
$3$-manifolds can be generated in the same way.
The basic ingredient of the algorithm is:
\begin{quote}
Start with some triangle and add further triangles
as long as no edge is contained in more than two
triangles. If this condition is violated, then backtrack.
A set of triangles is \emph{closed} if each of its edges
is contained in exactly two triangles. If the link of 
every vertex of a closed set of triangles is a circle,
then this set of triangles gives a triangulated surface: OUTPUT surface.
\end{quote}

From each equivalence class of combinatorially equivalent
triangulations (with respect to relabeling the
vertices) we list only the \emph{canonical} triangulation, 
the labeled triangulation which 
has the lexicographically smallest set of triangles in this class.
For every listed triangulation deg(1), the degree of vertex $1$,
must have minimum degree (since otherwise a lexicographically smaller
set of triangles can be obtained by relabeling the vertices) and the
triangulation must contain the triangles 
$$123,\, 124,\, 135,\, \dots,\, 1({\rm deg}(1)-1)({\rm deg}(1)+1),\, 1{\rm deg}(1)({\rm deg}(1)+1).$$

We enumerate the canonical triangulations in \emph{lexicographic order},
i.e., every listed triangulated surface is lexicographically smaller than the
next surface in the list.
With the objective to produce canonical triangulations we add the triangles 
during the backtracking in lexicographic order, that is, to the triangle $123$ 
we first add $124$ etc.\ to obtain a lexicographically ordered list of triangles.

We could wait until the list of triangles is a fully generated complex
before testing whether or not there are other combinatorially equivalent 
triangulations with lexicographically smaller lists of triangles.
However, we observe that, at each stage of adding triangles to obtain a 
canonical triangulation, the partial list of triangles is lexicographically
at least as small as any list obtained by relabeling the vertices.
We use this observation to prune the backtracking.

Whenever a new triangle is added to a partially generated complex,
we test whether the new complex can be relabeled to obtain a
lexicographically smaller labeling. If this is possible,
then the new partial complex will not lead to a canonical triangulation
and we backtrack.

\begin{table}
\small\centering
\defaultaddspace=0.1em
\caption{Backtracking steps in the case of $n=6$ vertices.}\label{tbl:back}
\begin{tabular*}{\linewidth}{@{}l@{\hspace{10mm}}l@{\extracolsep{5.75pt}}l@{}l@{}l@{}l@{}l@{}l@{\hspace{4.5mm}}l@{}l@{}}
\\\toprule
 \addlinespace
 \addlinespace
 \addlinespace
 \addlinespace
       & & \multicolumn{5}{@{}r}{Incomplete} & &  \\
Faces  & & \multicolumn{5}{@{}r}{vertices}   & & Reason for backtrack  \\
\midrule
 \addlinespace
 \addlinespace
 \addlinespace
 \addlinespace
 \addlinespace
123+124+134                              &&2&3&4&5&6&&   \\
123+124+134+234                          && & & &5&6&& surface complete \\
123+124+134                              &&2&3&4&5&6&&   \\
123+124+134+235                          &&2&3&4&5&6&&  \\
123+124+134+235+245                      && &3&4&5&6&&  \\
123+124+134+235+245+345                  && & & & &6&& surface complete \\
123+124+134+235+245                      && &3&4&5&6&& \\
123+124+134+235+245+346                  && &3&4&5&6&& \\
123+124+134+235+245+346+356              && & &4&5&6&& \\
123+124+134+235+245+346+356+456          && & & & & && \emph{surface complete}!\\
123+124+134+235                          &&2&3&4&5&6&& \\
123+124+134+235+246                      &&2&3&4&5&6&& \\
123+124+134+235+246+256                  && &3&4&5&6&& \\
123+124+134+235+246+256+345              && & &4&5&6&& relabeling is smaller\\
123+124+134+235+246+256                  && &3&4&5&6&& \\
123+124+134+235+246+256+346              && &3& &5&6&& relabeling is smaller\\
123+124+134                              &&2&3&4&5&6&& \\
123+124+135+145                          &&2&3&4&5&6&& \\
123+124+135+145+234                      && &3&4&5&6&& degree of 2 too small\\
123+124+135+145                          &&2&3&4&5&6&& \\
123+124+135+145+235                      &&2& &4&5&6&& degree of 3 too small\\
123+124+135+145                          &&2&3&4&5&6&& \\
123+124+135+145+236                      &&2&3&4&5&6&& \\
123+124+135+145+236+245                  &&2&3& &5&6&& degree of 4 too small\\
123+124+135+145+236                      &&2&3&4&5&6&& \\
123+124+135+145+236+246                  && &3&4&5&6&& \\
123+124+135+145+236+246+345              && &3&4& &6&& degree of 5 too small\\
123+124+135+145+236+246                  && &3&4&5&6&& \\
123+124+135+145+236+246+356              && & &4&5&6&& \\
123+124+135+145+236+246+356+456          && & & & & && \emph{surface complete}!\\
123+124+135+145                          &&2&3&4&5&6&& \\
123+124+135+146+156                      &&2&3&4&5&6&& \\
123+124+135+146+156+234                  && &3&4&5&6&& degree of 2 too small\\
123+124+135+146+156                      &&2&3&4&5&6&& \\
123+124+135+146+156+235                  &&2& &4&5&6&& degree of 3 too small\\
123+124+135+146+156                      &&2&3&4&5&6&& \\
123+124+135+146+156+236                  &&2&3&4&5&6&& \\
123+124+135+146+156+236+245              &&2&3&4&5&6&& \\
123+124+135+146+156+236+245+256          && &3&4&5&6&& \\
123+124+135+146+156+236+245+256+345      && &3&4& &6&& \\
123+124+135+146+156+236+245+256+345+346  && & & & & && \emph{surface complete}!\\
123+124+135+146+156+236                  &&2&3&4&5&6&& \\
123+124+135+146+156+236+246              && &3& &5&6&& degree of 2 too small\\
123+124+135+146+156                      &&2&3&4&5&6&& \\

 \addlinespace

 \addlinespace
 \addlinespace
 \addlinespace
 \addlinespace
 \bottomrule
\end{tabular*}
\end{table}

If there is a closed vertex $v$ such that $deg(v) < deg(1)$ then
there is a lexicographically smaller labeling; otherwise,
we search for such relabelings by
\begin{compactitem}
\item going through all closed vertices, $v$, 
      for which $\deg(v) = \deg(1)$,
\item and for each edge $vw$ we relabel $v$ as $1$ and $w$ as $2$,
\item thereafter we relabel the two vertices adjacent to the edge $vw$ 
      to be $3$ and $4$ (two choices).
\item Then we can extend the new labeling in a lexicographic smallest way.
\end{compactitem}

Table~\ref{tbl:back} displays the backtracking in case of $n=6$ vertices.
As a simplifying step we start not only with the triangle $123$,
but with the smallest possible completed vertex-star of size~$3$ 
of the vertex $1$, i.e., with the triangles $123+124+134$. 
The next smallest triangle is $234$ which closes the surface.
However, the resulting surface (the boundary of the tetrahedron)
has $4<6=n$ vertices and is therefore discarded.

Let $K$ be a partial complex and let $k$
be the smallest vertex of $K$ for which its vertex-star is
not closed. Since we add the new triangles in lexicographic order,
the next triangle to be added necessarily contains the vertex $k$.
In particular, the intersection of the new triangle with the  
current partial complex is not empty. Therefore,
\emph{every partial complex is connected}.

\begin{figure}
\begin{center}
\psfrag{1}{1}
\psfrag{2}{2}
\psfrag{3}{3}
\psfrag{4}{4}
\psfrag{5}{5}
\psfrag{6}{6}
\psfrag{7}{7}
\psfrag{8}{8}
\includegraphics[width=.40\linewidth]{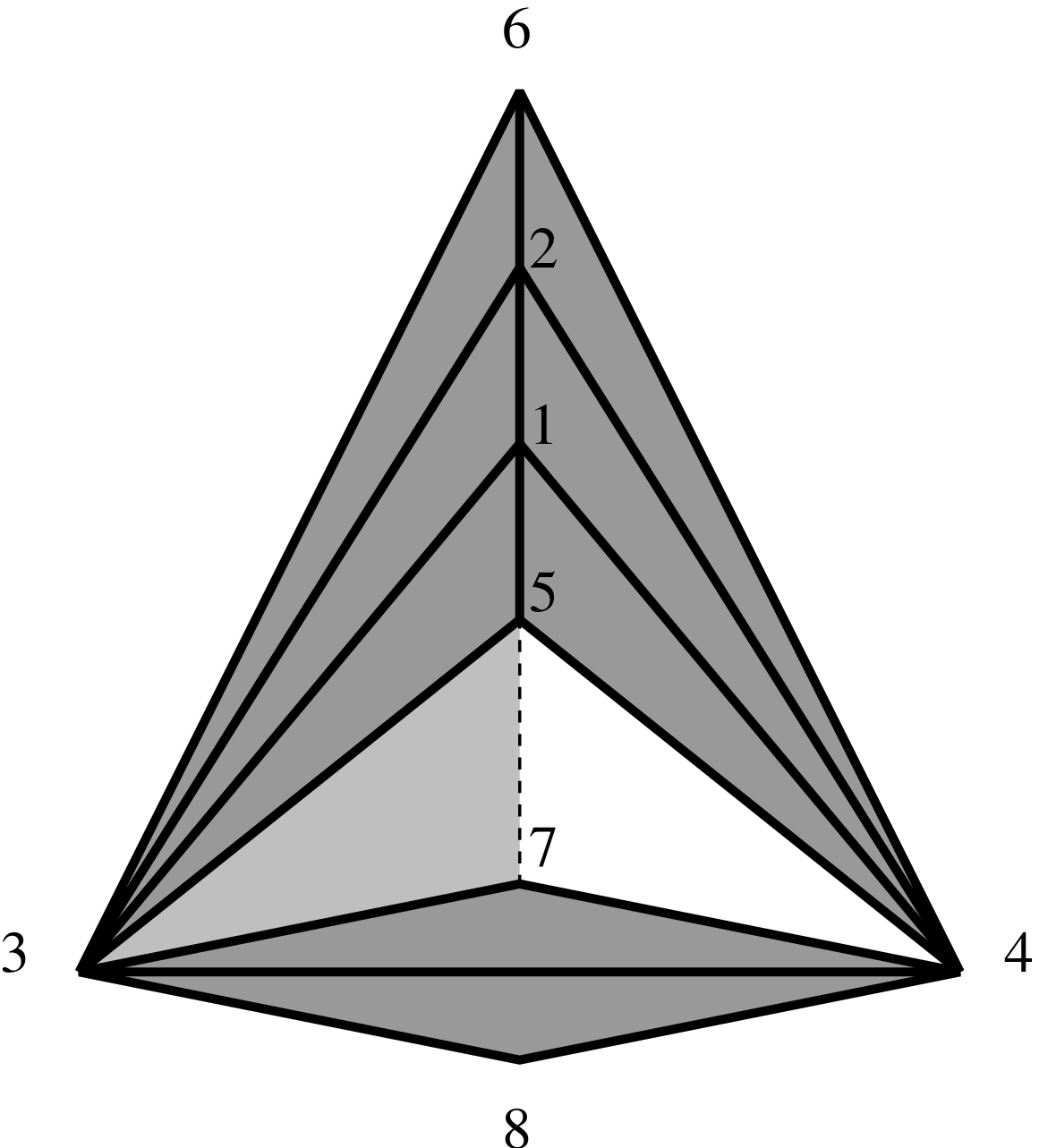}
\end{center}
\caption{Two strongly connected components that are joined by the triangle $357$.}
\label{fig:components}
\end{figure}

In a new triangle $klm$ at most the vertex $m$ has not yet been 
used as a vertex in the partial complex $K$.
(In the current vertex-star of $k$ in $K$ there
are at least two non-closed edges, say, $kr$ and $ks$.
Suppose $l$ and $m$ have not yet been used in $K$
and suppose $krx$ is the triangle that closes the edge $kr$.
Since $r,s<l,m$, it follows that $krx$ is lexicographically
smaller than $klm$ and is thus added to $K$ first, contradiction.) 
If $m$ is a new vertex
and $l$ is an existing vertex smaller than any
neighbor of $k$ on an unclosed edge, then $klm$ intersects
$K$ only in the vertices $k$ and $l$. Thus, $klm$ 
has no neighboring triangle in $K$. 
In other words, $klm$ (temporarily) forms a new strongly connected
component. For example, let the partial complex $K$
consist of the triangles
$$123,\,\,124,\,\,135,\,\,145,\,\,236,\,\,246,$$
to which the new component 
$$347,\,\,348$$
is added (see Figure~\ref{fig:components}), which then is connected to the first
component and closed to a triangulation of ${\mathbb R}{\bf P}^2$
by the triangles
$$357,\,\,368,\,\,458,\,\,467,\,\,567,\,\,568.$$
Thus, \emph{partial complexes are not necessarily strongly connected}.

\begin{prop}\label{prop:strongly}
Partial complexes are strongly connected upon completion of the link of a vertex.
\end{prop}

{\bf Proof.} Let $k$ be the smallest vertex for which its vertex-star is
not closed. We want to show that  at the time the last triangle is added 
to close the star of $k$ the resulting partial complex is strongly connected.
For this, we proceed by induction on $k$. 

First, we close the star of 1, 
which is a disc and therefore strongly connected.
We next assume that the partial complex which is obtained after closing the star of the vertex $k$
is strongly connected. Let $m > k$ be the next smallest vertex for
which its vertex-star is not yet closed. At the time we will have closed
the star of $m$, the star of $m$ is a disc (since otherwise we would
discard the respective partial complex). Since, by the induction
hypothesis, the partial complex after closing the star of $k$
was strongly connected, it follows that the partial complex after 
closing the star of $m$ is also strongly connected (because the star of $m$
contains at least one triangle that was present in the previous
partial complex).\\\mbox{}\hfill$\Box$

Let again $k$ be the smallest vertex for which its vertex-star is
not closed, $K$ be the current partial complex, and
$klm$ be the next triangle that is added to $K$.
In order for $klm$ to start a new strongly connected component,
$l$ has to be a vertex of the boundary of $K$ that is not (yet)
adjacent to $k$ and that is (by lexicographic minimality)
smaller than all other vertices of the boundary of $K$ 
to which $k$ is already adjacent.
The vertex $m$ might be a vertex of the boundary of $K$ or not.
If $m$ is a boundary vertex, then for $klm$ to start a new 
strongly connected component the vertices $k$, $l$, and $m$
are not pairwise adjacent in $K$. If $m$ is not a boundary vertex,
then we can choose (by lexicographic minimality) 
$m=|V(K)|+1$, with $V(K)$ the vertex set of $K$. 
The next triangle to be added to $K+kl(|V(K)|+1)$
is $kl(|V(K)|+2)$. The resulting strongly connected component $kl(|V(K)|+1)+kl(|V(K)|+2)$
cannot grow further. In the next step, either yet another
strongly connected component is started or the first strongly connected 
component is extended or joined to a later strongly connected component.

If during the enumeration of all $3$-manifolds with $n$ vertices 
the tetrahedron $klmr$ starts a new strongly connected component, 
then again $k$ is the currently smallest vertex for which its vertex-star is
not closed and $l$ belongs to the boundary of $K$, but is not
adjacent to $k$ in $K$. There are three cases for the vertices
$m$ and $r$. Either both belong to the boundary of $K$, in which case
$k$, $l$, $m$, and $r$ are not pairwise adjacent in $K$,
or only $m$ belongs to the boundary of $K$, in which case
$k$, $l$, and $m$ are not pairwise adjacent in $K$ and $r=|V(K)|+1$,
or both $m$ and $r$ do not belong to the boundary of $K$ and
$m=|V(K)|+1$ and $r=|V(K)|+2$.
Thus, the new strongly connected component consists
after its completion either of only the tetrahedron
$klmr$, of the two tetrahedra $klm(|V(K)|+1)$ and $klm(|V(K)|+2)$,
or of the join of the edge $kl$ with a circle that consists of $s$
edges $(|V(K)|+1)(|V(K)|+2)$, $(|V(K)|+1)(|V(K)|+3)$, $(|V(K)|+2)(|V(K)|+4)$,\dots, 
$(|V(K)|+s-1)(|V(K)|+s)$, where $|V(K)|+s\leq n$.

Proposition~\ref{prop:strongly} and the above analysis of
the strongly connected components explains why isomorphism-free
lexicographic enumeration is fast, but not as fast as
the generation of triangulations from irreducible triangulations.

In the latter approach one starts with the (finite) set of irreducible triangulations
of a surface from which triangulations with more vertices are obtained
by successive vertex-splitting. During this process the resulting
complexes are always proper triangulations of the initial surface. 
In other words, we stay within the class of triangulations of the surface.

In the lexicographic approach, the partial complexes do not necessarily
need to be strongly connected during the completion of the vertex-star
of the pivot vertex $k$. The possibility of more than one strongly
connected component leads to a  ``combinatorial explosion'' 
of the number of choices during the completion of the vertex
star $k$. Fortunately, the partial complexes become strongly connected
upon the completion of the vertex-star of~$k$. Thus the combinatorial explosion
happens locally, but not globally. Also, upon the completion of the
vertex-star of $k$ we can detect whether the link of $k$ is
indeed a triangulated $2$-sphere (or some other triangulated
$2$-manifold, in which case we discard the respective partial complex).

\bigskip

In the following sections we present our enumeration results
and corollaries thereof. In particular, we enumerated
all triangulated surfaces with 11 and 12 vertices
and all triangulated $3$-manifolds with $11$ vertices.

The algorithm was implemented as C programs which were executed on a cluster of
2GHz processors.  The total cpu time required was $20$ minutes to generate the
surfaces with $11$ vertices, $17$ days for the surfaces with $12$ vertices, 
and $170$ days for the $3$-manifolds with $11$ vertices.  
See \cite{Lutz_PAGE} for the program sources and lists of the examples.

\section{Triangulated Surfaces with 11 and 12 Vertices}
\label{sec:2d}

\begin{table}
\centering
\defaultaddspace=0.15em
\caption{Total numbers of triangulated surfaces with up to $12$ vertices.}\label{tbl:surf_total}
\begin{tabular}{@{}r@{\hspace{10mm}}r@{}}
\\\toprule
 \addlinespace
 \addlinespace
 \addlinespace
 \addlinespace
  $n$  &   Types \\
\midrule
\\[-4mm]
 \addlinespace
 \addlinespace
 \addlinespace
 \addlinespace
 
 \addlinespace
    4    &            1 \\
    5    &            1 \\
    6    &            3 \\
    7    &            9 \\
    8    &           43 \\
    9    &          655 \\
   10    &        42426 \\
   11    &     11590894 \\
   12    &  12561206794 \\
 \addlinespace

 \addlinespace
 \addlinespace
 \addlinespace
 \addlinespace
 \bottomrule
\end{tabular}
\end{table}

By Heawood's bound~\cite{Heawood1890}, at least 
$n\geq\Bigl\lceil\tfrac{1}{2}(7+\sqrt{49-24\chi (M)})\Bigl\rceil$ vertices
are needed to triangulate a (closed) surface of Euler characteristic $\chi (M)$. 
As shown by Ringel \cite{Ringel1955} and Jungerman and Ringel \cite{JungermanRingel1980}, 
this bound is tight, except in the cases of the orientable surface of genus~$2$, 
the Klein bottle, and the non-orientable surface of genus~$3$, 
for each of which an extra vertex has to be added.

Triangulations of surfaces with up to $8$ vertices were classified
by Datta \cite{Datta1999} and Datta and Nilakantan \cite{DattaNilakantan2002}.
By using (mixed) lexicographic enumeration, the second author
obtained all triangulations of surfaces with $9$ and $10$ vertices \cite{Lutz2005apre}.

We continued the enumeration with the isomorphism-free approach 
to lexicographic enumeration and were able to list all triangulated 
surfaces with $11$ and $12$ vertices. (Recently, Amendola \cite{Amendola2007pre}
independently generated all triangulated surfaces with $11$ vertices
by using genus-surfaces and isomorphism-free mixed-lexicographic enumeration.)

\begin{thm}\label{thm:surf}
There are precisely $11590894$ (combinatorially distinct) 
triangulated surfaces with $11$~vertices
and there are exactly $12561206794$ triangulated surfaces with $12$~vertices.
\end{thm}

The total numbers of triangulated surfaces with up to $12$ vertices
are given in Table~\ref{tbl:surf_total}. The numbers of triangulated
surfaces with $11$ and with $12$ vertices are listed in detail in the
Tables~\ref{tbl:surf_11} and \ref{tbl:surf_12}, respectively.

\begin{table}
\centering
\defaultaddspace=0.15em
\caption{Numbers of triangulated surfaces with $11$ vertices.}\label{tbl:surf_11}
\begin{tabular}{@{}r@{\hspace{8mm}}r@{\hspace{8mm}}r@{}}
\\\toprule
 \addlinespace
 \addlinespace
 \addlinespace
 \addlinespace
  Genus  &  Orientable  &  Non-orientable \\
\midrule
\\[-4mm]
 \addlinespace
 \addlinespace
 \addlinespace
 \addlinespace
 
 \addlinespace
    0    &         1249 &         -- \\
    1    &        37867 &      11719 \\
    2    &       113506 &      86968 \\
    3    &        65878 &     530278 \\
    4    &          821 &    1628504 \\
    5    &           -- &    3355250 \\
    6    &           -- &    3623421 \\
    7    &           -- &    1834160 \\
    8    &           -- &     295291 \\
    9    &           -- &       5982 \\
 \addlinespace

 \addlinespace
 \addlinespace
 \addlinespace
 \addlinespace
 \bottomrule
\end{tabular}
\end{table}

\begin{table}
\centering
\defaultaddspace=0.15em
\caption{Numbers of triangulated surfaces with $12$ vertices.}\label{tbl:surf_12}
\begin{tabular}{@{}r@{\hspace{8mm}}r@{\hspace{8mm}}r@{}}
\\\toprule
 \addlinespace
 \addlinespace
 \addlinespace
 \addlinespace
  Genus  &  Orientable  &  Non-orientable \\
\midrule
\\[-4mm]
 \addlinespace
 \addlinespace
 \addlinespace
 \addlinespace
 
 \addlinespace
    0    &        7595  &              -- \\
    1    &      605496  &          114478 \\
    2    &     7085444  &         1448516 \\
    3    &    25608643  &        16306649 \\
    4    &    14846522  &        99694693 \\
    5    &      751593  &       473864807 \\
    6    &          59  &      1479135833 \\
    7    &          --  &      3117091975 \\
    8    &          --  &      3935668832 \\
    9    &          --  &      2627619810 \\
   10    &          --  &       711868010 \\
   11    &          --  &        49305639 \\
   12    &          --  &          182200 \\
 \addlinespace

 \addlinespace
 \addlinespace
 \addlinespace
 \addlinespace
 \bottomrule
\end{tabular}
\end{table}

\begin{cor}
There are $821$ vertex-minimal triangulations of the orientable surface 
of genus $4$, and there are $295291$ and $5982$ vertex-minimal 
triangulations of the non-orientable surfaces of genus $8$ and $9$, 
respectively, with $11$ vertices.
\end{cor}

With a local search, Altshuler \cite{Altshuler1997} found
$59$ vertex-minimal \emph{neighborly} triangulations 
(i.e., with complete $1$-skeleton)
of the orientable surface of genus $6$ with $12$ vertices and
$40615$ neighborly triangulations with $12$ vertices 
of the non-orientable surface of genus $12$.
For the orientable surface of genus $6$ 
it was shown by Bokowski \cite{AltshulerBokowskiSchuchert1996}, \cite{Bokowski2006pre}
that Altshuler's list of $59$ vertex-minimal examples
is complete. For the non-orientable surface of genus $12$,
the $40615$ examples of Altshuler make up roughly one quarter 
of the exact number of $182200$ vertex-minimal
triangulations of this surface with $12$ vertices.

\begin{cor}
There are $751593$ vertex-minimal triangulations of the orientable surface 
of genus $5$, and there are $711868010$, $49305639$, and $182200$
vertex-minimal triangulations of the non-orientable surfaces of genus
$10$, $11$, and $12$, respectively, with $12$ vertices.
\end{cor}

The $182200$ vertex-minimal triangulations of the non-orientable surfaces 
of genus $12$ with $12$ vertices were previously generated 
by Ellingham and Stephens~\cite{EllinghamStephens2005}:
They used a modified isomorphism-free lexicographic enumeration
for the generation of all neighborly triangulations with $12$ and $13$
vertices. (There are 243088286 neighborly triangulations of the
non-orientable surface of genus $15$ with $13$ vertices \cite{EllinghamStephens2005}.)

\bigskip

Every $2$-dimensional simplicial complex (with $n$ vertices) is 
polyhedrally embeddable in~${\mathbb R}^5$, as it can be
realized as a subcomplex of the boundary complex of the cyclic
polytope $C(n,6)$; cf.\ Gr\"unbaum \cite[Ex.\ 25, p.\ 67]{Gruenbaum1967}.

However, not all triangulations of orientable surfaces are geometrically realizable 
in ${\mathbb R}^3$, i.e., with straight edges, flat triangles, and without self intersections:
Bokowski and Guedes de Oliveira \cite{BokowskiGuedes_de_Oliveira2000}
showed that one of the $59$ neighborly triangulations of the orientable
surface of genus $6$ is \emph{not} realizable in $3$-space. Recently,
Schewe \cite{Schewe2007} proved non-realizability in~${\mathbb R}^3$ 
for all the $59$ examples. Schewe further showed that for every orientable 
surface of genus $g\geq 5$ there are triangulations that cannot 
be realized in ${\mathbb R}^3$. 

Realizations for all vertex-minimal triangulations of the orientable surfaces
of genus $2$ and $3$ from \cite{Lutz2005apre} were obtained in \cite{Bokowski2006pre},
\cite{HougardyLutzZelke2006pre}, and \cite{Lutz2005apre},
and realizations of these triangulations with small coordinates in
\cite{HougardyLutzZelke2007a}, \cite{HougardyLutzZelke2007b}; 
see \cite{HougardyLutzZelke2006pre} for additional comments 
and further references on realizability.

The $821$ vertex-minimal triangulations of the orientable surface of genus $4$ 
from our enumeration were all found to be realizable \cite{HougardyLutzZelke2006pre}
as well as at least $15$ of the $751593$ vertex-minimal triangulations of the
orientable surface of genus $5$ with $12$ vertices. These results in
combination with the results of Schewe \cite{Schewe2007} led to:

\begin{conj} {\rm (Hougardy, Lutz, and Zelke, \cite{HougardyLutzZelke2006pre})}
Every triangulation of an orientable surface of genus
$1\leq g\leq 4$ is geometrically realizable.
\end{conj}

\section{Equivelar Surfaces}
\label{sec:equiv}

A particularly interesting class of triangulated surfaces are
\emph{equivelar simplicial maps}, i.e., triangulations
for which all vertices have the same vertex-degree $q$.
Equivelar simplicial maps are also called \emph{degree regular triangulations}
or \emph{equivelar triangulations}. 

In general, let a \emph{map} on a surface $M$ be a decomposition of $M$
into a finite cell complex and let $G$ be the $1$-skeleton of the map
on $M$. The graph $G$ of the map may have multiple edges, loops, vertices of degree $2$,
or even vertices of degree $1$; for example, the embedding of a tree
with $n$ vertices and $n-1$ edges on $S^2$ decomposes the $2$-sphere
into one polygon with $2n-2$ edges, which are identified pairwise. 
(Sometimes the graphs of maps are required to be connected finite
simple graphs, sometimes multiple edges are allowed but no loops,
and vertices are often required to have at least degree $3$;
see \cite{BrehmWills1993}, \cite{CoxeterMoser1957}, \cite{Ringel1974},~\cite{Vince2004}.)
A map is \emph{equivelar of type $\{p,q\}$} if $M$ is decomposed into
$p$-gons only with every vertex having degree $q$; cf.\
\cite{McMullenSchulzWills1982}, \cite{McMullenSchulzWills1983}. 
A map is \emph{polyhedral} if the intersection of any two of its polygons is either empty, a common vertex, 
or a common edge; see the surveys \ \cite{BrehmSchulte1997}, \cite{BrehmWills1993}.
An \emph{equivelar polyhedral map} is a map which is both equivelar and polyhedral. 

A map is \emph{regular} if it has a flag-transitive automorphism group. 
Regular maps therefore provide highly symmetric examples of equivelar maps;
see \cite{ConderDobcsanyi2001}, \cite{CoxeterMoser1957}, \cite{Wilson1976}.
Vertex-transitive maps and neighborly triangulations are further classes 
of equivelar surfaces that have intensively been studied in the
literature; cf.~\cite{Altshuler1997}, \cite{AltshulerBokowskiSchuchert1996}, 
\cite{EllinghamStephens2005}, \cite{JungermanRingel1980}, \cite{KoehlerLutz2005pre},
\cite{Ringel1955}.

Equivelar simplicial maps (as simplicial complexes) always are polyhedral.
By double counting of incidences between vertices and edges
as well as between edges and triangles, we have 
\begin{equation}
nq=2f_1=3f_2
\end{equation}
for equivelar triangulations, with $f_1$ and $f_2$ denoting the numbers of edges and $2$-faces, respectively.
By Euler's equation, we further have that
\begin{equation}\label{eq:equitri}
\chi(M)=n-f_1+f_2=n-\frac{nq}{2}+\frac{nq}{3}=\frac{n(6-q)}{6},
\end{equation}
or equivalently,
\begin{equation}
q=6-\frac{6\chi(M)}{n}.
\end{equation}
Since $q$ is a positive integer,
it follows that $n$ has to be a divisor of $6|\chi(M)|$ if $\chi(M)\neq 0$.
In particular, a surface $M$ of Euler characteristic $\chi(M)\neq 0$
has only finitely many equivelar triangulations. Table~\ref{tbl:equiv_val}
displays the possible values of $(n,q)$ for surfaces with $\chi(M)\geq -10$.

\begin{table}
\centering
\defaultaddspace=0.15em
\caption{Possible values of $(n,q)$ for equivelar triangulations with $\chi(M)\geq -10$.}\label{tbl:equiv_val}
\begin{tabular}{@{}r@{\hspace{8mm}}l@{}}
\\\toprule
 \addlinespace
 \addlinespace
 \addlinespace
 \addlinespace
  $\chi(M)$  &  $(n,q)$\\
\midrule
\\[-4mm]
 \addlinespace
 \addlinespace
 \addlinespace
 \addlinespace
 
 \addlinespace

    2  &  (4,3), (6,4), (12,5) \\
    1  &  (6,5) \\
    0  &  $(n,6)$, with $n\geq 7$ \\
  $-1$  &  -- \\
  $-2$  &                                    (12,7) \\
  $-3$  &                             (9,8), (18,7) \\
  $-4$  &                            (12,8), (24,7) \\
  $-5$  &                    (10,9), (15,8), (30,7) \\
  $-6$  &                    (12,9), (18,8), (36,7) \\
  $-7$  &                    (14,9), (21,8), (42,7) \\
  $-8$  &           (12,10), (16,9), (24,8), (48,7) \\
  $-9$  &                    (18,9), (27,8), (54,7) \\
 $-10$  &  (12,11), (15,10), (20,9), (30,8), (60,7) \\
 \addlinespace

 \addlinespace
 \addlinespace
 \addlinespace
 \addlinespace
 \bottomrule
\end{tabular}
\end{table}

In the case of neighborly triangulations we have $q=n-1$ and therefore $\chi(M)=\frac{n(7-n)}{6}$. 
It follows that $n\equiv 0,1,3,4\,{\rm mod}\,6$, where $n\geq 4$.
In the case $n=6k$ we have Euler characteristic $\chi(M)=-6k^2+7k$, and if $n=6k+1$ then $\chi(M)=-6k^2+5k+1$.
If $n=6k+3$ we have $\chi(M)=-6k^2+k+2$, and if $n=6k+4$ then $\chi(M)=-6k^2-k+2$.

Equivelar triangulations with up to $11$ vertices were classified 
by Datta and Nilakantan~\cite{DattaNilakantan2001}: there are $27$ such examples.
Datta and Upadhyay \cite{DattaUpadhyay2005} continued the classification 
of equivelar triangulations for the torus and the Klein bottle for up to $15$ vertices.
(Constructions of equivelar maps on the torus together with bounds on
their number were given in \cite{Altshuler1973};
for equivelar maps on the Klein bottle see \cite{Negami1984}.)
All equivelar polyhedral maps on the torus are vertex-transitive
\cite{BrehmKuehnel2006pre}, \cite{DattaUpadhyay2005}.
By isomorphism-free lexicographic enumeration, we obtained all equivelar triangulations 
of the torus and the Klein bottle for up to $100$ vertices; 
see Tables~\ref{tbl:equiv_torus} and \ref{tbl:equiv_klein}.

\begin{table}
\centering
\defaultaddspace=0.15em
\caption{Numbers of equivelar triangulations of the torus with up to $100$ vertices.}\label{tbl:equiv_torus}
\begin{tabular*}{\linewidth}{@{}c@{\extracolsep{8pt}}r@{\extracolsep{\fill}}r@{\extracolsep{\fill}}r@{\extracolsep{\fill}}r@{\extracolsep{\fill}}r@{\extracolsep{\fill}}r@{\extracolsep{\fill}}r@{\extracolsep{\fill}}r@{\extracolsep{\fill}}r@{\extracolsep{\fill}}r@{}}
\\\toprule
 \addlinespace
 \addlinespace
 \addlinespace
 \addlinespace
$k\backslash$Vertices  & {\small $10k+1$} & {\small $10k+2$} & {\small $10k+3$} & {\small $10k+4$} & {\small $10k+5$} & {\small $10k+6$} & {\small $10k+7$} & {\small $10k+8$} & {\small $10k+9$} & {\small $10(k+1)$}\\
\midrule
\\[-4mm]
 \addlinespace
 \addlinespace
 \addlinespace
 \addlinespace
 
 \addlinespace
     0   &  -- & -- & -- & -- & -- & -- &  1 &  1 &  2 &  1 \\
     1   &   1 &  4 &  2 &  2 &  4 &  5 &  2 &  5 &  3 &  6 \\ 
     2   &   6 &  4 &  3 & 11 &  5 &  5 &  7 &  9 &  4 & 11 \\
     3   &   5 & 11 &  8 &  7 &  8 & 16 &  6 &  8 & 10 & 16 \\ 
     4   &   6 & 15 &  7 & 13 & 14 & 10 &  7 & 24 & 10 & 14 \\  
     5   &  12 & 16 &  8 & 19 & 12 & 21 & 14 & 13 &  9 & 30 \\
     6   &  10 & 14 & 19 & 23 & 14 & 23 & 11 & 20 & 16 & 23 \\ 
     7   &  11 & 36 & 12 & 17 & 22 & 23 & 16 & 27 & 13 & 34 \\
     8   &  21 & 19 & 13 & 40 & 18 & 20 & 20 & 31 & 14 & 39 \\  
     9   &  20 & 27 & 22 & 22 & 20 & 47 & 16 & 27 & 27 & 37 \\
 \addlinespace

 \addlinespace
 \addlinespace
 \addlinespace
 \addlinespace
 \bottomrule
\end{tabular*}
\end{table}

\begin{table}
\centering
\defaultaddspace=0.15em
\caption{Numbers of equivelar triangulations of the Klein bottle with up to $100$ vertices.}\label{tbl:equiv_klein}
\begin{tabular*}{\linewidth}{@{}c@{\extracolsep{8pt}}r@{\extracolsep{\fill}}r@{\extracolsep{\fill}}r@{\extracolsep{\fill}}r@{\extracolsep{\fill}}r@{\extracolsep{\fill}}r@{\extracolsep{\fill}}r@{\extracolsep{\fill}}r@{\extracolsep{\fill}}r@{\extracolsep{\fill}}r@{}}
\\\toprule
 \addlinespace
 \addlinespace
 \addlinespace
 \addlinespace
$k\backslash$Vertices  & {\small $10k+1$} & {\small $10k+2$} & {\small $10k+3$} & {\small $10k+4$} & {\small $10k+5$} & {\small $10k+6$} & {\small $10k+7$} & {\small $10k+8$} & {\small $10k+9$} & {\small $10(k+1)$}\\
\midrule
\\[-4mm]
 \addlinespace
 \addlinespace
 \addlinespace
 \addlinespace
 
 \addlinespace
     0   & -- & -- & -- & -- & -- & -- & -- & -- &  1 &  1 \\
     1   & -- &  3 & -- &  1 &  3 &  2 & -- &  4 & -- &  4 \\ 
     2   &  3 &  1 & -- &  7 &  2 &  1 &  3 &  4 & -- &  8 \\ 
     3   & -- &  4 &  3 &  1 &  4 &  9 & -- &  1 &  3 &  8 \\ 
     4   & -- &  8 & -- &  4 &  7 &  1 & -- & 11 &  2 &  5 \\
     5   &  3 &  4 & -- &  8 &  4 &  8 &  3 &  1 & -- & 15 \\
     6   & -- &  1 &  7 &  6 &  4 &  8 & -- &  4 &  3 &  9 \\
     7   & -- & 15 & -- &  1 &  7 &  4 &  4 &  8 & -- & 12 \\
     8   &  5 &  1 & -- & 15 &  4 &  1 &  3 &  8 & -- & 16 \\
     9   &  4 &  4 &  3 &  1 &  4 & 15 & -- &  5 &  7 & 10 \\

 \addlinespace

 \addlinespace
 \addlinespace
 \addlinespace
 \addlinespace
 \bottomrule
\end{tabular*}
\end{table}

\begin{thm}
There are exactly $1357$ equivelar triangulations of the torus 
and $364$ equivelar triangulations of the Klein bottle 
with up to $100$ vertices, respectively. 
\end{thm}

Recently, Brehm and K\"uhnel \cite{BrehmKuehnel2006pre}
gave a detailed description of all equivelar triangulations 
of the torus. In particular, they obtained an explicit
formula for the number $T(n)$ of equivelar triangulations
with $n$ vertices (as well as for the number $Q(n)$
of equivelar polyhedral quadrangulations with $n$ vertices).


As observed by Datta and Upadhyay \cite{DattaUpadhyay2005},
there is an $n$-vertex equivelar triangulation of the Klein bottle
if and only if $n\geq 9$ is not prime.

Moreover, Datta and Upadhyay~\cite{DattaUpadhyay2006} determined that
there are exactly six equivelar triangulations of the orientable
surface of genus $2$ with $12$ vertices. As a consequence of Theorem~\ref{thm:surf}:
\begin{cor}
There are precisely $240914$ equivelar triangulations with $12$ vertices.
\end{cor}
Table~\ref{tbl:equiv} lists the numbers of simplicial equivelar maps with up to $12$ vertices.

\begin{table}
\centering
\defaultaddspace=0.15em
\caption{Numbers of simplicial equivelar maps with up to $12$ vertices.}\label{tbl:equiv}
\begin{tabular*}{\linewidth}{@{}c@{\extracolsep{16pt}}c@{\extracolsep{16pt}}r@{\extracolsep{16pt}}r@{\extracolsep{\fill}}c@{\extracolsep{16pt}}c@{\extracolsep{16pt}}r@{\extracolsep{16pt}}r@{}}
\\\toprule
 \addlinespace
 \addlinespace
 \addlinespace
 \addlinespace
  Vertices  &  Orient.  & Genus  &  Types   &  Vertices  &  Orient.  & Genus  &  Types \\
\midrule
\\[-4mm]
 \addlinespace
 \addlinespace
 \addlinespace
 \addlinespace
 
 \addlinespace
      4    &        $+$ &      0 &       1  &      12    &       $+$ &      0 &       1 \\
      6    &        $+$ &      0 &       1  &            &           &      1 &       4 \\
           &        $-$ &      1 &       1  &            &           &      2 &       6 \\
      7    &        $+$ &      1 &       1  &            &           &      3 &      34 \\
      8    &        $+$ &      1 &       1  &            &           &      4 &     112 \\
      9    &        $+$ &      1 &       2  &            &           &      5 &     103 \\
           &        $-$ &      2 &       1  &            &           &      6 &      59 \\
           &            &      5 &       2  &            &       $-$ &      2 &       3 \\
     10    &        $+$ &      1 &       1  &            &           &      4 &      28 \\
           &        $-$ &      2 &       1  &            &           &      6 &     500 \\
           &            &      7 &      14  &            &           &      8 &    9273 \\
     11    &        $+$ &      1 &       1  &            &           &     10 &   48591 \\
           &            &        &          &            &           &     12 &  182200 \\
 \addlinespace

 \addlinespace
 \addlinespace
 \addlinespace
 \addlinespace
 \bottomrule
\end{tabular*}
\end{table}

\bigskip

For an equivelar \emph{polyhedral} map of type $\{p,q\}$ the same computation as in Equation~(\ref{eq:equitri})
gives
\begin{equation}\label{eq:npq}
\chi(M)=n-f_1+f_2=n-\frac{nq}{2}+\frac{nq}{p}=nq\big(\frac{1}{p}+\frac{1}{q}-\frac{1}{2}\big).
\end{equation}
Thus, the sign of $\chi(M)$ is determined by the sign of $\frac{1}{p}+\frac{1}{q}-\frac{1}{2}$,
and vice versa. 

If $\frac{1}{p}+\frac{1}{q}-\frac{1}{2}>0$, then the only possible $\{p,q\}$-pairs 
are $\{3,3\}$, $\{3,4\}$, $\{3,5\}$, $\{4,3\}$, and $\{5,3\}$ for $S^2$ with $\chi(S^2)=2$,
with the boundaries of the tetrahedron, the octahedron, the icosahedron, the cube,
and the dodecahedron as the unique occurring examples, respectively, 
and $\{3,5\}$, $\{5,3\}$ for ${\mathbb R}{\bf P}^2$, with the vertex-minimal $6$-vertex
triangulation of ${\mathbb R}{\bf P}^2$ and its combinatorial dual as the only examples.

If $\frac{1}{p}+\frac{1}{q}-\frac{1}{2}=0$, then there are infinitely many
triangulations, quadrangulations, and hexangulations corresponding to
the pairs $\{3,6\}$, $\{4,4\}$, and $\{6,3\}$, respectively; 
see Brehm and K\"uhnel \cite{BrehmKuehnel2006pre} for more details.

In the case $\frac{1}{p}+\frac{1}{q}-\frac{1}{2}<0$
we write Equation~(\ref{eq:npq}) as
\begin{equation}\label{eq:q}
q=\frac{n-\chi(M)}{n}\cdot\frac{2p}{p-2},
\end{equation}
where $p$ and $q$ are positive integers greater or equal to $3$.

For a given surface $M$ of Euler characteristic $\chi(M)<0$ 
we next determine all triples $(p,q;n)$ which are admitted
by Equation~(\ref{eq:q}). Every equivelar polyhedral map has at least one $p$-gon
with $p$ vertices, i.e., we always have $n\geq p$. Furthermore, 
a vertex has $q$ distinct neighbors, which implies $n\geq q+1$.
The combinatorial dual of an equivelar polyhedral map of type $(p,q;n)$ 
is an equivelar polyhedral map of type $(q,p;\frac{nq}{p})$.
Moreover, in an equivelar polyhedral map of type $(p,q;n)$ 
the star of any vertex contains $q(p-3)+q+1=q(p-2)+1$ distinct vertices,
from which $n\geq q(p-2)+1\geq 4(p-2)+1=4p-7>2p$ follows for $q\geq p\geq 4$. 
In the case $p=3$ we have $q\geq 7$ for surfaces with $\chi(M)<0$ and therefore also $n>2p$.
If $q<p$, then for the dual maps of type $(q,p;\frac{nq}{p})$ we have
$\frac{nq}{p}>2q$, and thus again $n>2p$ for the maps of type $(p,q;n)$.

From Equation~(\ref{eq:q}) we see that $n$ is a divisor of $(n-\chi(M))2p$
and therefore a divisor of $2|\chi(M)|p$. Let $a$ be the gcd of $n$ and $p$, 
and let $k$ and $l$ be positive integers such that $n=ka$ and $p=la$. 
It follows from $n\,|\,2|\chi(M)|p$ that $k\,|\,2|\chi|$
and from $n>2p$ that $k>2l$. Thus $3\leq k\leq 2|\chi|$ and $1\leq l\leq \lfloor\frac{k-1}{2}\rfloor$,
that is, $k$ and $l$ can take only finitely many distinct values,
where $k\,|\,2|\chi|$ and ${\rm gcd}(k,l)=1$.
From 
\begin{equation}
q=\frac{ka-\chi(M)}{ka}\cdot\frac{2la}{la-2}=\frac{2la+\frac{2l|\chi(M)|}{k}}{la-2}=2+\frac{4+\frac{2l|\chi(M)|}{k}}{la-2},
\end{equation}
we see that there are only finitely many choices for $a$.
It follows, in particular, that for given $M$ with $\chi(M)<0$ 
there are only finitely many equivelar maps on $M$.

\bigskip

If $\chi(M)=-1$, then there are no admissible triples $(p,q;n)$.
Thus, there are no equivelar polyhedral maps on the non-orientable surface 
of genus $3$.

For $\chi(M)=-2$ the admissible triples are $(3,7;12)$ and $(7,3;28)$. 
Altogether, there are $12$ examples of equivelar polyhedral maps 
on the orientable surface of genus $2$ (the six simplicial examples 
from above and their simple duals); see \cite{DattaUpadhyay2006}.
\begin{cor}
There are exactly $56$ equivelar polyhedral maps on the non-orientable surface of genus~$4$,
$28$ of type $(3,7;12)$ and $28$ of type $(7,3;28)$.
\end{cor}
None of the examples of equivelar polyhedral maps with $\chi(M)=-2$
is regular.

For $\chi(M)=-3$ the admissible triples are
$(3,8;9)$, $(8,3;24)$, $(3,7;18)$, $(7,3;42)$, $(4,5;12)$, and $(5,4;15)$.
\begin{thm}
There are precisely $1403$ equivelar triangulations of the
non-orientable surface of genus $5$, two with $9$ vertices
and $1401$ with $18$ vertices.
\end{thm}
Furthermore, there are $4$ equivelar polyhedral maps on the non-orientable surface
of genus~$5$ of type $(4,5;12)$ \cite{LutzSulanke2007pre}.
One of these examples is regular; cf.\ \cite[p.~134]{Wilson1976}.

In the case $\chi(M)=-4$ we have the possibilities
$(3,8;12)$, $(8,3;32)$, $(3,7;24)$, $(7,3;56)$,
$(4,5;16)$, and $(5,4;20)$.

\begin{thm}
There are precisely $11301$ equivelar triangulations of the
orientable surface of genus $3$, $24$ with $12$ vertices
and $11277$ with $24$ vertices. Moreover, there are exactly
$601446$ equivelar triangulations of the
non-orientable surface of genus $6$, $500$ with $12$ vertices
and $600946$ with $24$ vertices.
\end{thm}
Exactly two of the equivelar triangulations with $\chi(M)=-4$
are regular,
Dyck's regular map (\cite{Dyck1880a}, \cite{Dyck1880b};
\cite{Bokowski1989}, \cite{BokowskiWills1988}, \cite{Brehm1987a}, 
\cite{SchulteWills1986a}, \cite{Sherk1959}, \cite{Wilson1976})
of type $(3,8;12)$ and Klein's regular map (\cite{Klein1879}; 
\cite{SchulteWills1985}, \cite{Sherk1959}, \cite{Wilson1976})
of type $(3,7;24)$.

There are $363$ equivelar polyhedral maps on the non-orientable
surface of genus $6$ of type $(4,5;16)$ \cite{LutzSulanke2007pre}
of which one is regular; cf.\ \cite[p.~139]{Wilson1976}.
Moreover, there are $43$ equivelar polyhedral maps on the orientable
surface of genus $3$ of type $(4,5;16)$, none of these are regular \cite{LutzSulanke2007pre}.

\bigskip

For neighborly triangulations of orientable surfaces the genus $g$
grows quadratically with the number of vertices $n$, i.e., $g=O(n^2)$.
However, the boundary of the tetrahedron and M\"obius' $7$-vertex torus \cite{Moebius1886}
are the only examples of neighborly triangulations of orientable
surfaces for which polyhedral realizations in ${\mathbb R}^3$
are known \cite{BokowskiEggert1991}, \cite{Csaszar1949}. 
In contrast, as mentioned above, all $59$ neighborly triangulations of the orientable
surface of genus $6$ with $12$ vertices are not realizable \cite{Schewe2007},
and it is expected that also all neighborly triangulations of
orientable surfaces with more vertices never are realizable.

McMullen, Schulz, and Wills~\cite{McMullenSchulzWills1983} constructed
polyhedral realizations in ${\mathbb R}^3$ of equivelar triangulations
of genus $g=O(n\log n)$, which, asymptotically, is the highest known 
genus $g(n)$ for geometric realizations of polyhedral maps. 
McMullen, Schulz, and Wills also gave infinite families 
of geometric realizations of equivelar polyhedral maps 
of the types $\{4,q\}$ and $\{p,4\}$. For further examples 
of geometric realizations of equivelar polyhedral maps 
of these types see
\cite{BokowskiWills1988},
\cite{McMullenSchulteWills1988},
\cite{McMullenSchulzWills1982}, 
\cite{SchulteWills1986b},
and \cite{Ziegler2004cpre}.

It is not known whether there are geometric realizations of equivelar
polyhedral maps of type $\{p,q\}$ for $p,q\geq 5$; cf.\ \cite{BrehmWills1993}. 
Examples of equivelar polyhedral maps of type $\{5,5\}$ and of type $\{6,6\}$ were first
given by Brehm~\cite{Brehm1990}. An infinite series of $\{k,k\}$-equivelar polyhedral maps 
was constructed by Datta~\cite{Datta2005}.

\section{Combinatorial 3-Manifolds with 11 Vertices}
\label{sec:3d}

\begin{table}
\small\centering
\defaultaddspace=0.15em
\caption{Combinatorial $3$-manifolds with up to $11$ vertices.}\label{tbl:ten3d_leq10}
\begin{tabular*}{\linewidth}{@{\extracolsep{\fill}}c@{\hspace{5mm}}r@{\hspace{5mm}}r@{\hspace{5mm}}r@{\hspace{5mm}}r@{\hspace{5mm}}r@{}}
\\\toprule
 \addlinespace
 \addlinespace
  Vertices$\backslash$Types & $S^3$ & $S^2\hbox{$\times\hspace{-1.62ex}\_\hspace{-.4ex}\_\hspace{.7ex}$}S^1$ & $S^2\!\times\!S^1$ & ${\mathbb R}{\bf P}^3$ &  All \\ 
\midrule
 \addlinespace
 \addlinespace
      5 &          1 &      --  &       -- &  -- &          1 \\
 \addlinespace
      6 &          2 &      --  &       -- &  -- &          2 \\
 \addlinespace
      7 &          5 &      --  &       -- &  -- &          5 \\
 \addlinespace
      8 &         39 &      --  &       -- &  -- &         39 \\
 \addlinespace
      9 &       1296 &        1 &       -- &  -- &       1297 \\
 \addlinespace
     10 &     247882 &      615 &      518 &  -- &     249015 \\
 \addlinespace
     11 &  166564303 &  3116818 &  2957499 &  30 &  172638650 \\ 
\bottomrule
\end{tabular*}
\end{table}

The boundary of the $4$-simplex triangulates the $3$-sphere
with $5$ vertices, and,
by work of Walkup \cite{Walkup1970}, the twisted sphere product 
$S^2\hbox{$\times\hspace{-1.62ex}\_\hspace{-.4ex}\_\hspace{.7ex}$}S^1$,
the sphere product $S^2\!\times\!S^1$, and the real projective $3$-space ${\mathbb R}{\bf P}^{\,3}$
can be triangulated vertex-minimally with $9$, $10$, and $11$ vertices,
respectively, while all other $3$-manifolds need at least $11$ vertices
for a triangulation.  By a result of Bagchi and Datta~\cite{BagchiDatta2005},
triangulations of ${\mathbb Z}_2$-homology spheres (different from $S^3$) 
require at least $12$ vertices.
In particular, at least $12$ vertices are needed to triangulate the lens space $L(3,1)$. 
A triangulation of $L(3,1)$ with this number of vertices was first found by Brehm~\cite{Brehm-pers}.
Otherwise, no bounds are known on the minimal numbers of vertices
of triangulated $3$-manifolds. 

Triangulations of $3$-manifolds with up to $10$ vertices were classified previously;
see \cite{Lutz2006apre} and the references given there.
With isomorphism-free lexicographic enumeration
we were able to obtain all triangulations with $11$ vertices.

\begin{thm}
There are precisely $172638650$ triangulated $3$-manifolds with $11$~vertices.
\end{thm}
Table~\ref{tbl:ten3d_leq10} lists the combinatorial and topological types
of the triangulations with up to~$11$ vertices. The numbers of triangulations
with $11$ vertices are displayed in detail in Table~\ref{tbl:ten3d_11}.

\begin{table}
\small\centering
\defaultaddspace=0.15em
\caption{Combinatorial $3$-manifolds with $11$ vertices.}\label{tbl:ten3d_11}
\begin{tabular*}{\linewidth}{@{\extracolsep{\fill}}c@{\hspace{5mm}}r@{\hspace{5mm}}r@{\hspace{5mm}}r@{\hspace{5mm}}r@{\hspace{5mm}}r@{}}
\\\toprule
 \addlinespace
 \addlinespace
$f$-vector$\backslash$Types &   $S^3$ & $S^2\hbox{$\times\hspace{-1.62ex}\_\hspace{-.4ex}\_\hspace{.7ex}$}S^1$ & $S^2\!\times\!S^1$ & ${\mathbb R}{\bf P}^3$ & All \\ 
\midrule
 \addlinespace
 \addlinespace
 (11,34,46,23) &        131 &          &          &     &        131 \\
 \addlinespace
 (11,35,48,24) &        859 &          &          &     &        859 \\
 \addlinespace
 (11,36,50,25) &       3435 &          &          &     &       3435 \\
 \addlinespace
 (11,37,52,26) &      11204 &          &          &     &      11204 \\
 \addlinespace
 (11,38,54,27) &      31868 &          &          &     &      31868 \\
 \addlinespace
 (11,39,56,28) &      82905 &          &          &     &      82905 \\
 \addlinespace
 (11,40,58,29) &     199303 &          &          &     &     199303 \\
 \addlinespace
 (11,41,60,30) &     447245 &          &          &     &     447245 \\
 \addlinespace
 (11,42,62,31) &     939989 &          &          &     &     939989 \\
 \addlinespace
 (11,43,64,32) &    1850501 &          &          &     &    1850501 \\
 \addlinespace
 (11,44,66,33) &    3413161 &      448 &      406 &     &    3414015 \\
 \addlinespace
 (11,45,68,34) &    5888842 &     3627 &     3521 &     &    5895990 \\
 \addlinespace
 (11,46,70,35) &    9463527 &    17065 &    16559 &     &    9497151 \\
 \addlinespace
 (11,47,72,36) &   14091095 &    54928 &    53839 &     &   14199862 \\
 \addlinespace
 (11,48,74,37) &   19288095 &   137795 &   134494 &     &   19560384 \\
 \addlinespace
 (11,49,76,38) &   23946497 &   278899 &   272671 &     &   24498067 \\
 \addlinespace
 (11,50,78,39) &   26344282 &   464328 &   451126 &     &   27259736 \\
 \addlinespace
 (11,51,80,40) &   24835145 &   626441 &   603950 &   1 &   26065537 \\
 \addlinespace
 (11,52,82,41) &   19130339 &   665845 &   630869 &   3 &   20427056 \\
 \addlinespace
 (11,53,84,42) &   11240196 &   525104 &   486378 &   6 &   12251684 \\
 \addlinespace
 (11,54,86,43) &    4457865 &   272672 &   244045 &   8 &    4974590 \\
 \addlinespace
 (11,55,88,44) &     897819 &    69666 &    59641 &  12 &    1027138 \\ \midrule
 \addlinespace
 \addlinespace
 Total:        &  166564303 &  3116818 &  2957499 &  30 &  172638650 \\ 
\bottomrule
\end{tabular*}
\end{table}

\begin{cor}
Let $M$ be a $3$-manifold different from $S^3$,  $S^2\hbox{$\times\hspace{-1.62ex}\_\hspace{-.4ex}\_\hspace{.7ex}$}S^1$,
$S^2\!\times\!S^1$, and ${\mathbb R}{\bf P}^3$  (which can be triangulated with $5$, $9$, $10$, and $11$ vertices,
respectively), then $M$ needs at least $12$ vertices for a triangulation. 
\end{cor}

\begin{cor}
There are exactly $30$ vertex-minimal triangulations of ${\mathbb R}{\bf P}^{\,3}$
with $11$ vertices. 
\end{cor}

\begin{cor}
Walkup's triangulation of ${\mathbb R}{\bf P}^{\,3}$  from \cite{Walkup1970}
is the unique vertex- and facet-minimal triangulation of ${\mathbb R}{\bf P}^{\,3}$
with $f=(11,51,80,40)$.
\end{cor}

\begin{cor}
The minimal number of vertices for triangulations of the orientable connected sum 
$(S^2\!\times\!S^1)\# (S^2\!\times\!S^1)$ and of the non-orientable connected sum
$(S^2\hbox{$\times\hspace{-1.62ex}\_\hspace{-.4ex}\_\hspace{.7ex}$}S^1)\# (S^2\hbox{$\times\hspace{-1.62ex}\_\hspace{-.4ex}\_\hspace{.7ex}$}S^1)$
is $12$.
\end{cor}
Examples of triangulations of the latter two manifolds with $12$ vertices
are given in \cite{Lutz2005bpre}. It is conjectured in
\cite{Lutz2005bpre} that for other $3$-manifolds, different from
the mentioned six examples, at least $13$ vertices are necessary
for a triangulation.

In \cite{Lutz2006apre}, all triangulated $3$-spheres with up to $10$ vertices
and all resulting simplicial $3$-balls with $9$ vertices were examined
with respect to shellability. The respective $3$-spheres all
turned out to be shellable, whereas $29$ vertex-minimal examples
of non-shellable $3$-balls were discovered with $9$ vertices;
see also \cite{Lutz2004a}.

\begin{cor}
All triangulated $3$-spheres with $11$ vertices are shellable.
\end{cor}

The smallest known example of a non-shellable $3$-sphere
has $13$ vertices \cite{Lutz2004c}. We believe that 
there are no non-shellable $3$-spheres with $12$ vertices.

\begin{cor}
There are $1831363502$ triangulated $3$-balls with $10$ vertices
of which $277479$ are non-shellable.
\end{cor}

For all triangulated $3$-spheres with up to $9$ vertices
and all neighborly $3$-spheres with $10$ vertices
a classification into polytopal and non-polytopal examples
was carried out mainly by Altshuler, Bokowski, and Steinberg;
see \cite{Lutz2005apre} for a survey and references.

\begin{prob}
Classify all simplicial $3$-spheres with $10$ and $11$ vertices
into polytopal and non-poly\-topal spheres.
\end{prob}

\bibliography{.}

\bigskip
\bigskip
\medskip

\small

\noindent
Thom Sulanke\\
Department of Physics\\
Indiana University\\
Bloomington\\
Indiana 47405\\
USA\\
{\tt tsulanke@indiana.edu}\\[2mm]

\medskip

\noindent
Frank H. Lutz\\
Technische Universit\"at Berlin\\
Fakult\"at II - Mathematik und Naturwissenschaften\\
Institut f\"ur Mathematik, Sekr.\ MA 3-2\\
Stra\ss e des 17.\ Juni 136\\
10623 Berlin\\
Germany\\
{\tt lutz@math.tu-berlin.de}

\end{document}